\documentclass[10pt,twoside]{classe-ma_f}
\usepackage{amssymb,amsbsy,amsmath,amsfonts,amssymb,amscd}
\usepackage[english,francais]{babel}
\usepackage{times}
\newtheorem{theoreme}{Th\'eor\`eme}[section]

\newtheorem{proposition}[theoreme]{Proposition}
\newtheorem{corollaire}[theoreme]{Corollaire}

\newtheorem{question}{\it Question}

\AuteurCourant{Alireza Abdollahi}
 \TitreCourant{Une condition nilpotence} 
\begin{document}
\selectlanguage{francais}
\title{Groupes satisfaisant une condition nilpotence}
\author{Alireza Abdollahi}
\address{\labelsep=2mm\leftskip=-5mm
   D\'epartement de Math\'ematiques, Universit\'e d'Ispahan,
 Ispahan 81746-73441, Iran.\\ Courriel:  a.abdollahi@sci.ui.ac.ir }
\maketitle
\begin{Resume}{%
Dans cette Note, nous \'etudions les groupes $G$ satisfont la
condition $(\mathcal{N},3)$, c'est-\`a-dire, si  chaque ensemble
de $n+1$ \'el\'ements de $G$ contient une paire $\{x,y\}$ telle
que le sous-groupe $\left(x,y\right>$ soit nilpotent. }
\end{Resume}
\selectlanguage{english}
\begin{Etitle}{Groups satisfying a nilpotence condition}\end{Etitle}
\begin{Abstract}{%
In this Note we study the groups $G$ satisfying condition
$(\mathcal{N},n)$, that is, every subset of $G$ with $n+1$
elements contains a pair $\{x,y\}$ such that the subgroup
$\left<x,y\right>$ is nilpotent.
 }\end{Abstract}


\par\medskip\centerline{\rule{2cm}{0.2mm}}\medskip
\setcounter{section}{0} \selectlanguage{francais}
\section{Introductions et r\'esultats}
Soient $\mathcal{X}$ une classe de groupes  et $n$ un entier
positif.  Nous dirons  qu'un groupe $G$ satisfait la condition
$(\mathcal{X},n)$ si  chaque ensemble de $n+1$ \'el\'ements de $G$
contient une paire $\{x,y\}$ telle que le sous-groupe
$\left(x,y\right>$ soit un $\mathcal{X}$-groupe. Par exemple, si
$\mathcal{X}$ est frem\'ee par passage aux sous-groupes, un groupe
qui est la r\'eunion de $n$ sous-groupes de $\mathcal{X}$
satisfait la condition $(\mathcal{X},n)$. Soit $\mathcal{N}$ la
classe des groupes nilpotents. Dans \cite{endi}, Endimioni a
prouv\'e qu'un groupe fini $G$ satisfaisant la condition
$(\mathcal{N},3)$ est nilpotent  et le groupe sym\'etrique $S_3$
satisfait la condition $(\mathcal{N},4)$. Une question qui se
pose naturellement est
\begin{question}
Y a-t-il un groupe satisfaisant la condition $(\mathcal{N},3)$ qui
n'est pas faiblement nilpotent? (on dit qu'un groupe est
faiblement nilpotent s'il satisfait la condition
$(\mathcal{N},1)$).
\end{question}
 Nous
donnons quelques r\'eponses partielles \`a la question ci-dessus. 
On sais que dans un groupe $G$ les
conditions suivantes sont \'equivalentes:\\
1) ~ $G$ est nilpotent.\\
2) $G$ est faiblement nilpotent.\\
3) ~ $G$ est le produit direct de ses sous-groupes de Sylow.\\
Il est bien connus que les conditions ci-dessus ne sont pas
\'equivalentes dans un groupe quel conque. Dans cette Note, on
preuve que si l'un groupe $G$ satisfaisant la condition
$(\mathcal{N},3)$, alors l'ensemble d'\'el\'ements de $G$ d'ordre
fini  est un sous-groupe de $G$ qui est le produit direct des
sous-groupes de Sylow de $G$ et le produit direct des
$2'$-sous-groupes de Sylow de $G$ est faiblement nilpotent (voir
la proposition \ref{p1}). Nous g\'en\'eralisons en outre le
r\'esultat d'Endimioni \cite{endi} \`a quelques grandes classes
des groupes contenant la classe des groupes finis
 comme des groupes
satisfaisant la condition maximale sur les sous-groupes ab\'eliens
et groupes  r\'esiduellement  ($2'$-fini) (voir le corollaire \ref{c1}).\\
Nous noterons $\mathcal{N}$ ou $\mathcal{N}_{\infty}$ la classe
des groupes nilpotents et si $\lambda$ un entier positif, nous
noterons $\mathcal{N}_{\lambda}$ la class des groupes nilpotents
de  classe au plus \'egale \`a $\lambda$. Notons qu'un groupe fini
satisfaisant la condition $(\mathcal{N},n)$, en fait, satisfait
la condition $(\mathcal{N}_{c},n)$ pour un entier positif $c$.
Ainsi une autre question qui se pose naturellement est celle pour
laquelle  nombre entier positif $c$, chaque groupe satisfaisant
la condition $(\mathcal{N}_c,3)$  est faiblement nilpotent. Dans
la dernier partie de cette Note, on \'etudie des groupes
satisfaisant $(\mathcal{N}_c,3)$ pour petit $c$.\\

Cette recherche a \'et\'e financ\'ee par le  projet (num\'ero
801031) de l'universit\'e d'Ispahan.

\section{D\'emonstrations}

\begin{proposition} \label{p1}
 Soient $G$ un groupe
satisfaisant la condition $(\mathcal{N}_{\lambda},3)$ o\`u
$\lambda$ un entier positif ou $\infty$ et $p,q$ deux nombres
premiers distincts. Alors
\begin{enumerate}
\item  $\left<x^2,y\right>\in \mathcal{N}_{\lambda}$ et
$\left<x,x^y\right>\in \mathcal{N}_{\lambda}$ pour tous $x,y\in
G$. En particulier, si $\lambda<\infty$ alors $G$ est un groupe de
$(\lambda+1)$-Engel et si $\lambda=\infty$ alors $G$ est un
groupe d'Engel.
\item Si $x$ est un $p$-\'el\'ement et $y$ est un $q$-\'el\'ement
de $G$, alors $xy=yx$.
\item  l'ensemble de $p$-\'el\'ements  de $G$ est un sous-groupe
$T_p$. En particulier, l'ensemble des \'el\'ements d'ordres finis
est un sous-groupe $T$ de $G$ qui est le produit direct des $T_p$
. \item   $S=\text{Dr}_{p\not=2} T_p$ satisfait la condition
$(\mathcal{N}_{\lambda},1)$.
\end{enumerate}
\end{proposition}
\begin{demo} (1)  Consid\'erons  l'ensemble
$\{x,y,xy,x^2y\}$. Il est ais\'e de voir que  chaque sous-groupe
engendr\'e par deux \'el\'ements de l'ensemble ci-dessus est
\'egal \`a $\left<x^2,y\right>$ ou \`a $\left<x,y\right>$. Alors
par hypoth\`ese, $\left<x,y\right>\in\mathcal{N}_{\lambda}$ ou
$\left<x^2,y\right>\in\mathcal{N}_{\lambda}$. Mais
$\left<x^2,y\right>\leq \left<x,y\right>$, donc certainement
$\left<x^2,y\right>\in\mathcal{N}_{\lambda}$.\\ Pour prouver
$\left<x,x^y\right>\in\mathcal{N}_{\lambda}$, consid\'erons
l'ensemble $\{x,y,xy,yx\}$ et encore notons que chaque sous-groupe
engendr\'e par deux \'el\'ements de l'ensemble ci-dessus est
\'egal \`a $\left<x,y\right>$ ou \`a $\left<xy,xy\right>$. Alors
par hypoth\`ese, $\left<x,y\right>\in\mathcal{N}_{\lambda}$ ou
$\left<xy,xy\right>\in\mathcal{N}_{\lambda}$. Mais
$\left<xy,yx\right>\leq \left<x,y\right>$, donc certainement
$\left<xy,yx\right>\in\mathcal{N}_{\lambda}$. Maintenant,
supposons que $a, b\in G$, $x=b$ et $y=b^{-1}a$; alors on a
$\left<a,a^b\right>\in\mathcal{N}_{\lambda}$.\\
 Notons que si $\left<x,x^y\right>$ est nilpotent de classe au plus $c$
alors  $[(x^y)^{-1},_{c} x]=1$, donc
$1=[[y,x],_{c}x]=[y,_{c+1}x]$,  ce qui ach\`eve la
d\'emonstration de la partie (1).\\

(2)  Nous avons $$x^{-1}x^y=[x,y]=y^{-x}y.$$ D'apr\`es la partie
(1),  les deux sous-groupes $\left<x^y,x\right>$ et
$\left<y^x,y\right>$
 sont nilpotents. Pour cela, $[x,y]$ est non seulement un $p$-\'el\'ement
mais aussi un $q$-\'el\'ement.
 Il en r\'esulte que $[x,y]=1$.

 (3)  Supposons que $x,y$ sont deux $p$-\'el\'ements de $G$. Si $p\not
=2$, alors $\left<x^2,y\right>=\left<x,y\right>$. Alors, d'apr\`es
la partie (1), $\left<x,y\right>$ est nilpotent. Ainsi $xy$ est un
$p$-\'el\'ement de $G$. Maintenant, supposons que $p=2$ et
$x^{2^n}=1$, pour un certain $n\in\mathbb{N}$. Ainsi nous avons
$$(xy)^{2^n}=x^{-2^n}(xy)^{2^n}=y^{x^{2^n-1}}y^{x^{2^n-2}}\cdots
y^xy.$$ Puisque $\left<a,a^b\right>$ est nilpotent pour tous
$a,b\in G$, d'apr\`es la partie (1),
 on a $a^ba$ est un $2$-\'el\'ement, lorsque $a$ est un $2$-\'el\'ement de
$G$. Alors $y^xy$ est
 un $2$-\'el\'ement et donc aussi $y^{x^3}y^{x^2}y^xy=(y^xy)^{x^2}(y^xy)$
est un $2$-\'el\'ement.
 Par r\'ecurrence, nous obtenons que
 $(xy)^{2^n}=y^{x^{2^n-1}}y^{x^{2^n-2}}\cdots y^xy$
 est un $2$-\'el\'ement et alors  $xy$ est un $2$-\'el\'ement.\\
 Alors, $T_p$ est un sous-groupe de $G$. Maintenant supposons $a,b$ deux
\'el\'ements d'ordres finis dans
 $G$, $a=\prod_{p}a_p$  et $b=\prod_{p}b_p$  sont leurs decompositions dans
un produit de $p$-\'el\'ements,
 alors d'apr\`es la partie (2),  $$ab=\prod_{p}a_pb_p$$
 et les facteurs  commutent deux \`a deux. Alors $ab$ est d'ordre fini.
Il en r\'esulte que $T$ est un sous-groupe et $T$ est le produit
direct des $T_p$.\\

(4)  Supposons que $x,y\in S$; nous devons montrer que
$\left<x,y\right>\in\mathcal{N}_{\lambda}$. Mais l'ordre de $x$
est premier \`a $2$, donc $\left<x^2,y\right>=\left<x,y\right>$.
Alors d'apr\`es la partie (1),
$\left<x,y\right>\in\mathcal{N}_{\lambda}$.
\end{demo}

\begin{corollaire}\label{c1}
Soit $G$ un groupe satisfaisant la condition
$(\mathcal{N}_{\lambda},3)$ o\`u $\lambda$ est un entier positif
ou $\infty$.
\begin{enumerate}
\item Si $G$ satisfait la condition maximal sour les sous-groupes
ab\'eliens, alors $G$ est nilpotent. \item Si $G$ est
r\'esiduaillement  ($2'$-groupe fini), alors $G$ satisfait la
condition $(\mathcal{N}_{\lambda},1)$.
\item Si $\lambda$ est fini et $G$ est r\'esiduaillement fini alors $G$ localement nilpotent.
\end{enumerate}
\end{corollaire}
\begin{demo}
  (1) ~ D'apr\`es la proposition \ref{p1} (1), $G$ est un groupe d'Engel.
 Mais d'apr\`es un r\'esultat de Peng (Voir le th\'eor\`eme 7.21 de \cite{Rob2}), $G$ est
 nilpotent.\\
 (2) ~ Soit $x,y$ deux \'el\'ements de $G$. Alors
 $K:=\left<x,y\right>$ est r\'esiduaillement ($2'$-groupe fini). D'apr\`es la
 proposition \ref{p1} (1), $\left<x^2,y\right>$ est nilpotent (de classe au plus \'egale \`a
 $\lambda$).
 Soit $c$ la classe de $\left<x^2,y\right>$ ($c\leq\lambda$, si $\lambda<\infty$).
 Si $N\trianglelefteq
 K$ tel que $K/N$ est un $2'$-groupe fini, donc
 $\left<x^2,y\right>/N=K/N$. Maintenant $K/N$ est un groupe
 nilpotent de classe au plus \'egale \`a $c$. Alors $\gamma_{c+1}(K)\leq N$ et
 puisque $K$ est r\'esiduaillement ($2'$-groupe fini), on a
 $\gamma_{c+1}(K)=1$. Il implique que $G$ satisfait la condition
 $(\mathcal{N}_{\lambda},1)$.\\
 (3) ~ D'apr\`es la proposition \ref{p1}, $G$ est un  groupe de
 $(\lambda+1)$-Engel. Alors d'apr\`es un r\'esultat de Wilson
 \cite{W}, $G$ est localement nilpotent.
\end{demo}

\begin{proposition}
Soit $G$ un groupe satisfaisant la condition $(\mathcal{N}_n,3)$.
\begin{enumerate}
\item Si $n=1$ alors $G$ est nilpotent de classe au plus \'egale \`a
$2$.
\item Si $n=2$ alors $G$ est un groupe $3$-Engel. En particulier
$G$ est groupe  r\'esoluble de classe au plus $3$ et localement
nilpotent.
\item Si $n=3$ alors $G$ est un groupe $4$-Engel  localement nilpotent.
\end{enumerate}
\end{proposition}
\begin{demo}
(1)~ D'apr\`es la proposition \ref{p1} (1), $G^2\leq Z(G)$. Mais
$G'\leq G^2$, il en r\'esult que $G\in\mathcal{N}_2$.\\
(2)~ D'apr\`es la proposition \ref{p1} (2),
  $G$ est 3-Engel. Maintenant un r\'esultat bien connu de Heineken
(voir \cite{H}) implique que $G$ est localement nilpotent. Mais
d'apr\`es la proposition \ref{p1}, $G^2\leq R_2(G)$, o\`u
$$R_2(G)=\{x\in G \;|\; [x,y,y]=1 \;\;\forall y\in G \}$$  est l'ensemble les \'el\'ements 2-Engel \`a droite de
$G$. D'apr\`es un r\'esultat de Kappe (voir le corollaire 1 dans
la page 44 de \cite{Rob2}) $R_2(G)$ est un sous-groupe normal de
$G$. Aussi d'apr\`es un r\'esultat de Levi (voir le corollaire 2
dans la page 45 de \cite{Rob2}) on a $R_2(G)\in\mathcal{N}_3$.
Alors $G^2$ est r\'esoluble de classe au plus 2. Il en r\'esult
que $G$ est r\'esoluble de classe au plus 3. \\
(3)~ D'apr\`es la proposition \ref{p1},   $G$ est un groupe
$4$-Engel. Aussi d'apr\`es la proposition \ref{p1},
$\left<x^2,y\right>\in\mathcal{N}_3$ pour chaque $x,y\in G$. Il
en r\'esult que $[y,[x^2,[x^2,y]]]=1$. Alors le lemme 8 de
\cite{H} implique que $x^2\in HP(G)$ (le radical de Hirsch-Plotkin
de $G$). Il en r\'esult que $G^2$ est localement nilpotent. Donc
$G$ est un groupe (localement nilpotent)-par ab\'elien.
Maintenant un r\'esultat de Plotkin \cite{Plot} qui dit chaque
groupe de radical et Engel est localement nilpotent ach\`eve la
d\'emonstration.
\end{demo}


\end{document}